\documentclass[12pt,reqno]{article}
\usepackage{amssymb,amscd,amsmath,amsthm,color}
\newtheorem{theorem}{Theorem}[section]
\newtheorem{lemma}[theorem]{Lemma}
\newtheorem{cor}[theorem]{Corollary}
\newtheorem{prop}[theorem]{Proposition}

\usepackage{enumerate,amssymb}
\theoremstyle{definition}
\newtheorem{definition}[theorem]{Definition}

\theoremstyle{remark}

\numberwithin{equation}{section}

\def\bF{\mathbb{F}}

\def\bK{\mathbb{K}}
\def\bM{\mathbb{M}}

\def\bR{\mathbb{R}}
\def\diag{\mathrm{diag}}
\def\bJ{\mathbb{J}}
\def\bR{\mathbb{R}}
\def\bB{\mathbb{B}}

\begin{document}
\baselineskip=15pt

\title{Matrix inequalities from a two variables functional}

\author{ Jean-Christophe Bourin and  Eun-Young Lee{\footnote{
Research  supported by Basic Science Research Program through the National Research Foundation of Korea(NRF)
funded by the Ministry of Education, Science and Technology (2013-R1A1A2059872).}}  }

\date{ }

\maketitle

\vskip 10pt\noindent
{\small
{\bf Abstract.} We introduce a  two variables norm functional and establish its joint log-convexity. This entails  and improves many remarkable matrix inequalities, most of them related to the  log-majorization theorem of Araki. In particular: {\it   if $A$ is a  positive semidefinite matrix and $N$ is a normal matrix, $p\ge 1$ and $\Phi$ is a sub-unital positive linear map, then $|A\Phi(N)A|^p$ is weakly log-majorized by $A^p\Phi(|N|^p)A^p$}. This far extension of Araki's theorem (when $\Phi$ is the identity and $N$ is positive) complements some recent results of  Hiai and contains several special interesting cases such as a triangle inequality for normal operators and some extensions of the Golden-Thompson trace inequality.  Some applications to Schur products are also obtained. 
}

\section{Log-majorization and log-convexity }

\vskip 5pt\noindent
  Matrices   are regarded as  non-commutative extensions of
 scalars and functions.  Since matrices do not commute in general, most scalars identities cannot be brought to the matrix setting, however they sometimes have a matrix version, which is not longer an identity but an inequality. These kind of inequalities are of fundamental importance in our understanding of the noncommutative world of matrices. A famous, fifty years old example of such an inequality is the  Golden-Thompson trace inequality:  for Hermitian  $n$-by-$n$  matrices $S$ and $T$,
$$
{\mathrm{Tr\,}} e^{S+T}  \le {\mathrm{Tr\,}} e^{S/2} e^T e^{S/2} .
$$
A decade after,  Lieb and Thirring \cite{LT} obtained a stronger, remarkable trace inequality: for all positive semidefinite $n$-by-$n$ matrices $A,B\in\bM_n^+$ and all integers $p\ge 1$,
\begin{equation}\label{LT}
{\mathrm{Tr\,}} (ABA)^p \le {\mathrm{Tr\,}} A^pB^pA^p.
\end{equation}
This was  finally extended some fifteen years later  by Araki \cite{Araki} as a very important theorem in matrix analysis and its applications. Given $X,Y\in \bM_n^+ $, we write $X\prec_{w\!\log} Y$ when the series of $n$ inequalities holds, 
$$
\prod_{j=1}^k \lambda_j(X) \le \prod_{j=1}^k \lambda_j(Y)
$$
for  $k=1,\ldots n$, where $\lambda_j(\cdot)$ stands for the eigenvalues arranged in decreasing order. If further equality occurs for $k=n$, we  write $X\prec_{\log} Y$. Araki's  theorem  considerably strenghtens the Lieb-thirring trace inequality as the beautiful  log-majorization
\begin{equation}\label{araki}
 (ABA)^p \prec_{\log} A^pB^pA^p
\end{equation}
for all real numbers $p\ge 1$. In particular, this ensures \eqref{LT} for all $p\ge 1$.

Log- and weak log-majorization relations play a fundamental role in matrix analysis, a basic one for normal operators $X,Y\in\bM_n$ asserts that
\begin{equation}\label{folklore}
 |X+Y| \prec_{w\!\log} |X|+|Y|.
\end{equation}
This useful version of the triangle inequality belongs to the folklore and is a byproduct of Horn's inequalities, see the proof of \cite[Corollary 1.4]{bourin-IJM}.

This article  aims to provide new matrix inequalities containing \eqref{araki} and \eqref{folklore}. These inequalities are given in Section 2. The first part  dealing with positive operators is closely related to a recent paper of Hiai \cite{Hiai2}. The second part of Section 2 considers normal operators and contains our main theorem (Theorem \ref{cornormal}), mentioned in the Abstract.

Our main idea, and technical tool, is Theorem \ref{th1} below. It establishes the log-convexity  of a two variables functional. Fixing one variable in this functional yields a generalization of $\eqref{araki}$ involving a third matrix $Z\in\bM_n $, of the form
\begin{equation*}
 (AZ^*BZA)^p \prec_{w\!\log} A^pZ^*B^pZA^p.
\end{equation*}
We will also derive the following  weak log-majorization which contains both \eqref{araki} and \eqref{folklore} and thus unifies these two inequalities.

\vskip 10pt
\begin{prop}\label{th0}  Let $A\in\bM_n^+$ and  let $X,Y\in\bM_n$  be normal. Then, for all $p\ge 1,$
\begin{equation*}
 |A(X+Y)A|^p \prec_{w\!\log} 2^{p-1}A^p(|X|^p+|Y|^p)A^p.
\end{equation*}
\end{prop}

\vskip 5pt
Letting $X=Y=B$ in Proposition \ref{th0} we have \eqref{araki}, more generally,
\begin{equation}\label{araki-normal}
|AXA|^p \prec_{\log} A^p|X|^pA^p
\end{equation}
for all $A\in\bM_n^+$ and normal matrices $X\in\bM_n$. When $X$ is Hermitian, this  was noted by Audenaert \cite[Proposition 3]{Audenaert}.
If $A$ is the identity and $p=1$,  Proposition \ref{th0} gives \eqref{folklore}.  From \eqref{araki-normal}    follows several nice inequalities for the matrix exponential, due to  Cohen and al.\ \cite{Cohen1}, \cite{Cohen2}, including the  Golden-Thompson trace inequality and the elegant relation
\begin{equation}\label{exponential}
|e^Z| \prec_{\log} e^{{\mathrm{Re}}\,Z}
\end{equation}
for all matrices $Z\in\bM_n$, where ${\mathrm{Re}}\,Z=(Z+Z^*)/2$, \cite[Theorem 2]{Cohen1}.

Fixing the other variable in Theorem \ref{th1} below entails a H\"older inequality due to Kosaki. Several 
matrix versions of an inequality of Littlewood related to H\"older's inequality will be also obtained. 

The two variables in Theorem \ref{th1} are essential and  reflect  a construction with the perspective of a convex function. Recall that a  norm  on   $\bM_n$ is symmetric  whenever  $\| UAV\|=\| A\|$ for all $A\in\bM_n$ and all unitary  $U,V\in\bM_n$.  For $X,Y\in\bM_n^+$, the condition $X\prec_{w\!\log} Y$ implies $\| X\| \le
\| Y\|$ for all symmetric norms. We  state our log-convexity theorem.

\vskip 10pt
\begin{theorem}\label{th1} Let $A,B\in\bM_n^+$ and $Z\in\bM_n$. Then, for all symmetric norms and   $\alpha>0$, the map
$$
(p,t) \mapsto \left\| \left|A^{t/p}ZB^{t/p}\right|^{\alpha p} \right\| 
$$
is jointly log-convex on $(0,\infty)\times(-\infty,\infty)$.
\end{theorem}

\vskip 5pt
Here, if $A\in\bM_n^+$ is not invertible, we  naturally define for  $t\ge 0$,  $A^{-t}:=(A+F)^{-t}E$ where $F$ is the projection onto the nullspace of $A$ and $E$ is the range projection of $A$.

The next two sections present many hidden consequences of Theorem \ref{th1}, several of them extending \eqref{araki} and/or \eqref{folklore}, for instance,
$$
\left|A\frac{T+T^*}{2}A\right|^p \prec_{w\!\log} A^p\frac{|T|^p+|T^*|^p}{2}A^p
$$
for all $A\in\bM_n^+$, $p\ge 1$, and any $T\in\bM_n$. The proof of Theorem \ref{th1} is in Section 4. The  last section provides a version of Theorem \ref{th1}   for  operators acting on an infinite dimensional Hilbert space.

\section{Araki type inequalities}

\subsection{With positive operators}

\vskip 5pt To obtain new Araki's type inequalities, we  fix $t=1$ in Theorem \ref{th1} and thus use  the following special case.

\vskip 5pt
\begin{cor}\label{cornew} Let $A,B\in\bM_n^+$ and $Z\in\bM_n$. Then, for all symmetric norms and   $\alpha>0$, the map
$$
p \mapsto \left\| \left|A^{1/p}ZB^{1/p}\right|^{\alpha p } \right\| 
$$
is  log-convex on $(0,\infty)$.
\end{cor}

\vskip 5pt
We may now state a series of corollaries extending Araki's theorem.

\vskip 5pt
\begin{cor}\label{cor1}  Let $A,B\in\bM_n^+$ and $p\ge 1$. Let $Z\in\bM_n$ be a contraction. Then, for all symmetric norms and   $\alpha>0$, 
$$
  \| (AZ^*BZA)^{\alpha p} \| \le \| (A^{p}Z^*B^{p}ZA^{p})^{\alpha}\|.
$$
\end{cor}

\vskip 5pt Let $I$ be the identity of $\bM_n$. A matrix $Z$ is contractive, or a contraction, if $Z^*Z\le I$, equivalently if its operator norm satisfies $\| Z\|_{\infty}\le 1$.

\vskip 5pt
\begin{proof} The function $f(p)=\| | B^{1/p}ZA^{1/p} |^{2\alpha p} \|$ is log-convex, hence convex on $(0,\infty)$, and bounded since $Z$ is contractive, $0\le f(p) \le \| B\|_{\infty}^{2\alpha}
\| A\|_{\infty}^{2\alpha} \| I\|$. Thus $f(p)$ is nonincreasing, so $f(1)\ge f(p)$ for all $p\ge1$. Replacing $B$ by $B^{p/2}$ and $A$ by $A^{p}$ completes the proof.
\end{proof}

\vskip 5pt
Let $\|\cdot\|_{\{k\}}$, $k=1,\ldots, n$, denote the normalized Ky Fan $k$-norms on $\bM_n$,
$$
\| T\|_{\{k\}} =\frac{1}{k}\sum_{j=1}^k \lambda_j(|T|).
$$
Since, for all $A\in\bM_n^+$,
$$
\lim_{\alpha\to 0^+} \| A^{\alpha}\|_{\{k\}}^{1/\alpha} =\left\{\prod_{j=1}^k \lambda_j(A)\right\}^{1/k}
$$
we obtain from Corollary \ref{cor1} applied to the normalized Ky Fan $k$-norms, with $\alpha\to 0^+$,  a striking weak-log-majorization extending Araki's theorem.

\vskip 5pt
\begin{cor}\label{corstriking}  Let $A,B\in\bM_n^+$ and $p\ge 1$. Then,  for all contractions $Z\in \bM_n$,
$$
   (AZ^*BZA)^p \prec_{w\!\log}   A^{p}Z^*B^{p}ZA^{p}.
$$
\end{cor}

\vskip 5pt
If $Z=I$, we have the determinant equality and thus Araki's log-majorization \eqref{araki}.
Corollary \ref{cor1} and \ref{corstriking} are equivalent. Our proof of these extensions of Araki's theorem follows from  the two variables technic of Theorem \ref{th1}.  It's worth mentioning that Fumio Hiai  also obtained Corollary \ref{corstriking} in the beautiful note \cite{Hiai2}. Hiai's approach is based on some subtle estimates for the operator geometric mean.

For $X,Y\in \bM_n^+ $, the notation $X\prec^{w\!\log} Y$ indicates that the series of $n$ inequalities holds, 
$$
\prod_{j=1}^k \nu_j(X) \ge \prod_{j=1}^k \nu_j(Y)
$$
for  $k=1,\ldots n$, where $\nu_j(\cdot)$ stands for the eigenvalues arranged in increasing order. The following so-called super weak-log-majorization is another extension of Araki's theorem.
A matrix  $Z$ is expansive when  $Z^*Z\ge I$.

\vskip 5pt
\begin{cor}\label{corsuper}  Let $A,B\in\bM_n^+$ and $p\ge 1$. Then,  for all expansive matrices $Z\in \bM_n$,
$$
   (AZ^*BZA)^p \prec^{w\!\log}   A^{p}Z^*B^{p}ZA^{p}.
$$
\end{cor}

\vskip 5pt
\begin{proof} By a limit argument, we may assume invertibility of $A$ and $B$. Taking inverses, and using that $Z^{-1}$ is contractive, Corollary \ref{corsuper} is then equivalent to Corollary \ref{corstriking}.
\end{proof}

\vskip 5pt Corollaries \ref{corstriking}-\ref{corsuper}  imply a host of trace inequalities. We say that a continuous function $h:[0,\infty)\to (-\infty,\infty)$ is e-convex, (resp.\ e-concave), if $h(e^t)$ is convex, (resp.\ concave) on $(-\infty,\infty)$.  For instance, for all $\alpha >0$, $t\mapsto \log(1+t^{\alpha})$ is e-convex, while $t\mapsto \log(t^{\alpha}/(t+1))$ is e-concave. 
The equivalence between Corollary \ref{corstriking} and Corollary \ref{cortrace} below  is a basic property of majorization discussed in any monograph on this topic such as \cite{Bhatia} and \cite{HP}.

\vskip 25pt
\begin{cor}\label{cortrace}  Let $A,B\in\bM_n^+$,  $Z\in\bM_n$, and $p\ge 1$. 
\begin{itemize}
\item[{\rm(a)}] If $Z$ is contractive and $f(t)$ is e-convex and nondecreasing, then
$$
  {\mathrm{Tr\,}} f((AZ^*BZA)^p) \le  {\mathrm{Tr\,}} f( A^{p}Z^*B^{p}ZA^{p}).
$$
\item[{\rm(b)}]  If $Z$ is expansive and $g(t)$ is e-concave and nondecreasing, then
 $$
  {\mathrm{Tr\,}} g((AZ^*BZA)^p) \ge  {\mathrm{Tr\,}} g( A^{p}Z^*B^{p}ZA^{p}).
$$
\end{itemize}
\end{cor}

\vskip 5pt 
We will propose in Section 4 a proof of Theorem \ref{th1} making use of antisymmetric tensor powers, likewise in the proof of Araki's log-majorization. We will also indicate another, more elementary way, without antisymmetric tensors. The antisymmetric tensor technic goes back to Hermann Weyl, cf.\ \cite{Bhatia}, \cite{HP}. We  use it  to derive our next corollary.

\vskip 5pt
\begin{cor}\label{corlim}  Let $A,B\in\bM_n^+$ and  $Z\in\bM_n$. For each $j=1,\ldots, n$, the function defined on $(0,\infty)$
$$
p\mapsto \lambda_j^{1/p}( A^{p}Z^*B^{p}ZA^{p})
$$ 
converges  as $p\to\infty$.
\end{cor}

\vskip 5pt
\begin{proof} We may assume that $Z$ is contractive. As in the proof of Corollary \ref{cor1} we then see that the function
$
g(p)=\lambda_1^p( A^{1/p}Z^*B^{1/p}ZA^{1/p})
$
is log-convex and bounded, hence nonincreasing on $(0,\infty)$. Therefore $g(p)$ converges as $p\to 0$ and so $g(1/p)$ converges as $p\to \infty$. Thus 
$
p\mapsto  \lambda_1^{1/p}( A^{p}Z^*B^{p}ZA^{p})
$
converges
 as $p\to\infty$. Considering $k$-th antisymmetric tensor products, $k=1,\ldots, n$, we infer the convergence of
$$
p\mapsto\prod_{j=1}^k \lambda_j^{1/p}( A^{p}Z^*B^{p}ZA^{p})=\lambda_1^{1/p}\left((\wedge^kA)^p\wedge^kZ^*
(\wedge^kB)^p\wedge^kZ(\wedge^kA)^p\right)
$$
 and so, the convergence of
$
p\mapsto \lambda_j^{1/p}( A^{p}Z^*B^{p}ZA^{p})
$
as $p\to\infty$, for each $j=1,2,\ldots.$
\end{proof}

\vskip 5pt
 When $Z=I$, Audenaert and Hiai \cite{Aud-Hiai} recently gave a remarkable improvement of Corollary \ref{corlim} by showing that $p\mapsto (A^pB^pA^p)^{1/p}$ converges in $\bM_n$ as $p\to \infty$. We do not know whether such a reciprocal Lie-Trotter limit still holds with a third matrix $Z$ as in Corollary \ref{corlim}.

It is possible to state Corollary \ref{corstriking} in a  stronger form involving a positive linear map $\Phi$. Such a map is called sub-unital when $\Phi(I)\le I$.

\vskip 5pt
\begin{cor}\label{corstronger}  Let $A,B\in\bM_n^+$ and $p\ge 1$. Then,  for all positive linear, sub-unital map $\Phi: \bM_n\to\bM_n$,
$$
   (A\Phi(B)A)^p \prec_{w\!\log}   A^{p}\Phi(B^{p})A^{p}.
$$
\end{cor}

\vskip 5pt
\begin{proof} We may assume (the details are given, for a more general class of maps,  in the proof of Corollary 
\ref{poslin}) that
$$
\Phi(X) =\sum_{i=1}^m Z_i^* XZ_i
$$
where $m=n^2$ and $Z_i\in\bM_n$, $i=1,\ldots, m$, satisfy $\sum_{i=1}^m Z^*_iZ_i \le I$. Corollary \ref{corstronger} then follows from Corollary \ref{corstriking} applied to the  operators $\tilde{A},\tilde{B},\tilde{Z}\in \bM_{mn}$,
\begin{equation*}\label{eqblock}
\tilde{A}=\begin{pmatrix} A& 0_n&\cdots &0_n \\
0_n &0_n&\cdots &0_n \\
\vdots &\vdots &\ddots &\vdots \\
0_n& 0_n&\cdots &0_n \\
\end{pmatrix}, \
\tilde{B}=\begin{pmatrix} B& 0_n&\cdots &0_n \\
0_n &B&\cdots &0_n \\
\vdots &\vdots &\ddots &\vdots \\
0_n& 0_n&\cdots &B \\
\end{pmatrix}, \
\tilde{Z}=\begin{pmatrix} Z_1& 0_n&\cdots &0_n \\
Z_2 &0_n&\cdots &0_n \\
\vdots &\vdots &\ddots &\vdots \\
Z_m& 0_n&\cdots &0_n \\
\end{pmatrix}
\end{equation*}
where $0_n$ stands for the zero matrix in $\bM_n$.
\end{proof}

\vskip 5pt
  Corollary \ref{corstronger} can be applied for the Schur product $\circ$ (i.e., entrywise product) in $\bM_n$.

\vskip 5pt
\begin{cor}\label{corschur}  Let $A,B,C\in\bM_n^+$ and $p\ge 1$. If $C$ has all its diagonal entries less than or equal to one, then \
$$
   (A(C\circ B)A)^p \prec_{w\!\log}   A^{p}(C\circ B^{p})A^{p}.
$$
\end{cor}

\vskip 5pt
\begin{proof} The map $X\mapsto C\circ X$ is a positive linear, sub-unital map on $\bM_n$.
\end{proof}

\vskip 5pt
Corollary \ref{corschur} with the matrix $C$ whose entries are all equal to one is Araki's log-majorization. With $C=I$, Corollary \ref{corschur} is already an interesting extension of 
Araki's theorem as we may assume that $B$ is diagonal in \eqref{araki}. We warn the reader that the super weak-log-majorization, for  $A,B\in\bM_n^+$ and $p\ge 1$,
$
   (A(I\circ B)A)^p \prec^{w\!\log}   A^{p}(I\circ B^{p})A^{p}
$
does not hold, in fact, in general, $\det^ 2 I\circ B< \det I\circ B^2 $. 

\subsection{With normal operators}

To obtain Proposition \ref{th0} we need the following generalization of Corollary \ref{corstronger}.

\vskip 5pt
\begin{theorem}\label{cornormal}  Let $A\in \bM_n^+$ and let $N\in\bM_m$ be normal.   Then,  for all positive linear, sub-unital maps $\Phi: \bM_m\to\bM_n$, and $p\ge 1$,
$$
   |A\Phi(N)A|^p \prec_{w\!\log}   A^{p}\Phi(|N|^{p})A^{p}.
$$
\end{theorem}

\vskip 5pt
\begin{proof} By completing, if necessary, our matrices $A$ and $N$ with some 0-entries, we may assume that $m=n$ and then, as in the proof of Corollary \ref{corstronger}, that $\Phi$ is a congruence map with a contraction $\tilde{Z}$, $\Phi(X)=\tilde{Z}X\tilde{Z}^*$. Now, we have with the polar decomposition $N=U|N|$,
\begin{align*}
|A\tilde{Z}N\tilde{Z}^*A| &= \left|A\tilde{Z}|N|^{1/2}U|N|^{1/2}\tilde{Z}^*A\right| \\
&\prec_{\log} A\tilde{Z}|N|\tilde{Z}^*A
\end{align*}
by using  Horn's log-majorization $|XKX^*|\prec_{w\!\log} XX^*$ for all $X\in\bM_n$ and all contractions $K\in\bM_n$. Hence, from Corollary \ref{corstriking}, for all $p\ge 1$,
$$
|A\tilde{Z}N\tilde{Z}^*A|^p \prec_{\log} \left|A\tilde{Z}|N|\tilde{Z}^*A\right|^p  \prec_{w\!\log} A^p\tilde{Z}|N|^p\tilde{Z}^*A^p
$$
which completes the proof.
\end{proof}

\vskip 10pt
We are in a position to prove Proposition \ref{th0}  whose $m$-variables version is given here.

\vskip 10pt
\begin{cor}\label{th0a}  Let $A\in\bM_n^+$ and  let $X_1,\cdots,X_m\in\bM_n$  be normal. Then, for all $p\ge 1,$
\begin{equation*}
 \left|A\left(\sum_{k=1}^m X_k\right)A\right|^p \prec_{w\!\log} m^{p-1}A^p\left(\sum_{k=1}^m |X_k|^p\right)A^p.
\end{equation*}
\end{cor}

\vskip 10pt\noindent
\begin{proof}  Applying Theorem \ref{cornormal} to
$N=X_1\oplus \cdots\oplus X_m$
and to the unital, positive linear map $\Phi:\bM_{mn}\to \bM_n$, 
$$
\begin{pmatrix} S_{1,1}&\cdots&S_{1,m} \\ \vdots& \ddots &\vdots \\ S_{m,1}&\cdots&S_{m,m}\end{pmatrix}\mapsto \frac{1}{m}\sum_{k=1}^m S_{k,k}
$$
yields 
\begin{equation*}
 \left|A\frac{\sum_{k=1}^m X_k}{m}A\right|^p \prec_{w\!\log} A^p\frac{\sum_{k=1}^m |X_k|^p}{m}A^p
\end{equation*}
which is equivalent to the desired inequality.
\end{proof}

\vskip 5pt
A special case of Corollary \ref{th0a} deals with the Cartesian decomposition of an arbitrary matrix.

\vskip 5pt
\begin{cor}\label{corcartesian} Let $X,Y\in \bM_n$ be Hermitian. Then, for all $p\ge 1$, 
\begin{equation*}
 \left|A(X+iY)A\right|^p  \prec_{w\!\log} 2^{p-1}A^p(|X|^p+|Y|^p)A^p
\end{equation*}
where the constant $2^{p-1}$ is the best possible.
\end{cor}

\vskip 5pt
To check that $2^{p-1}$ is optimal, take $A=I\in\bM_{2n}$ and pick any two-nilpotent matrix,
$$
X+iY=\begin{pmatrix} 0&T \\ 0&0\end{pmatrix}.
$$

For a single normal operator, Corollary \ref{th0a} gives \eqref{araki-normal} as we have equality for the determinant. This entails the following remarkable log-majorization for the matrix exponential.

\vskip 5pt
\begin{cor}\label{corthompson} Let $A,B\in \bM_n$. Then,  
\begin{equation*}
\left| e^{A+B}  \right|\prec_{\log} e^{{\mathrm{Re}\,} A/2}  e^{{\mathrm{Re}\,}B}
e^{{\mathrm{Re}\,} A/2}.
\end{equation*}
\end{cor}

\vskip 5pt
 Corollary \ref{corthompson} contains \eqref{exponential} and shows that when $A$ and $B$ are Hermitian we have the famous Thompson log-majorization, \cite[Lemma 6]{Thompson},
\begin{equation}\label{gtlog}
 e^{A+B}  \prec_{\log} e^{ A/2}  e^{B}
e^{ A/2}
\end{equation}
which entails
\begin{equation*}
 \| e^{A+B} \| \le \| e^{ A/2}  e^{B}
e^{ A/2}\|
\end{equation*}
for all symmetric norms. For the operator norm this is Segal's inequality while for the trace norm this is the Golden-Thompson inequality. Taking the logarithms in \eqref{gtlog}, we have a classical majorization between $A+B$ and $\log e^{A/2}e^Be^{A/2}$. Since $t\mapsto |t|$ is convex, we infer, replacing $B$ by $-B$ that
\begin{equation*}
 \| A-B \| \le \| \log (e^{ A/2}  e^{-B}
e^{ A/2})\|
\end{equation*}
for all symmetric norms.
For the Hilbert-Schmidt norm, this is the Exponential Metric Increasing inequality, reflecting the nonpositive curvature of the positive definite cone  with its  Riemannian structure (\cite[Chapter 6]{Bhatia2}).

 Corollary \ref{corthompson} follows from \eqref{araki-normal} combined with the Lie Product Formula \cite[p.\ 254]{Bhatia} as shown in the next proof. Note that  Corollary \ref{corthompson} also follows  from Cohen's log-majorization \eqref{exponential} combined  with Thompson's log-majorization \eqref{gtlog}, thus we do not pretend to originality.

\vskip 5pt
\begin{proof} We have a Hermitian matrix $C$ such that, using the Lie Product Formula,
$$
e^{A+B}=e^{{\mathrm{Re}\,}A +{\mathrm{Re}\,} B +iC}=\lim_{n\to+\infty} \left( e^{({\mathrm{Re}\,}A+{\mathrm{Re}\,}B)/2n}e^{iC/n}e^{({\mathrm{Re}\,}A
+{\mathrm{Re}\,}B)/2n}\right)^{n}.
$$
On the other hand, by \eqref{araki-normal}, for all $n\ge 1$,
$$
\left| \left( e^{({\mathrm{Re}\,}A+{\mathrm{Re}\,}B)/2n}e^{iC/n}e^{({\mathrm{Re}\,}A
+{\mathrm{Re}\,}B)/2n}\right)^{n} \right| \prec_{\log} e^{{\mathrm{Re}\,}A +{\mathrm{Re}\,} B}
$$
so that
$$
\left|e^{A+B}\right| \prec_{\log}  e^{{\mathrm{Re}\,}A +{\mathrm{Re}\,} B}.
$$
Using again the Lie Product Formula,
$$
 e^{{\mathrm{Re}\,}A +{\mathrm{Re}\,} B} =\lim_{n\to+\infty} \left( e^{{\mathrm{Re}\,}A/2n} e^{{\mathrm{Re}\,}B/n}e^{{\mathrm{Re}\,}A/2n}
\right)^{n},
$$
combined with \eqref{araki-normal} (or \eqref{araki}) completes the proof.
\end{proof}

\vskip 5pt
 Theorem \ref{cornormal} is the main result of Section 2 as all the other results in this section are special cases.
One more elegant extension of Araki's inequality  follows, involving an arbitrary matrix.

\vskip 5pt
\begin{cor}\label{corlast} Let $A\in \bM_n^+$ and $p\ge 1$. Then, for any $T\in\bM_n$,
\begin{equation*}
 \left|A\frac{T+T^*}{2}A\right|^p \prec_{w\!\log} A^p\frac{|T|^p+|T^*|^p}{2}A^p.
\end{equation*}
\end{cor}

\vskip 5pt
\begin{proof} It suffices to apply  Theorem \ref{cornormal} to 
$$
N=\begin{pmatrix} 0& T \\ T^*&0\end{pmatrix}
$$
and to the unital, positive linear map $\Phi:\bM_{2n}\to\bM_n$,
$$
\begin{pmatrix} B&C  \\ D&E \end{pmatrix}\mapsto \frac{B+C+D+E}{2}.
$$
\end{proof}

\vskip 5pt
We apply Theorem \ref{cornormal} to Schur products in the next two corollaries.

\vskip 5pt
\begin{cor}\label{corschurlast}  Let $A\in\bM_n^+$ and  let $X,Y\in\bM_n$  be normal. Then, for all $p\ge 1,$
\begin{equation*}
 |A(X\circ Y)A|^p \prec_{w\!\log} A^p(|X|^p\circ |Y|^p)A^p.
\end{equation*}
\end{cor}

\vskip 5pt
\begin{proof} We need to see the Schur product as a positive linear map,
$$
X\circ Y= \Phi(X\otimes Y)
$$
where $\Phi:\bM_n\otimes\bM_n\to \bM_n$ merely consists in extracting a principal submatrix. Setting $N=X\otimes Y$ in Theorem \ref{cornormal} completes the proof.
\end{proof}

\vskip 5pt
We  note that Corollary \ref{corschurlast} extends \eqref{araki-normal} (with $X$ in diagonal form  and  $Y=I$) and contains the classical log-majorization for normal operators,
$$|X\circ Y|\prec_{w\\log} |X|\circ |Y|.$$
As a last illustration of the scope of Theorem \ref{cornormal} we have the following result.

\vskip 5pt
\begin{cor}\label{corverylast} Let $A\in \bM_n^+$ and $p\ge 1$. Then, for any $T\in\bM_n$,
\begin{equation*}
 \left|A(T\circ T^*)A\right|^p \prec_{w\!\log} A^p(|T|^p\circ|T^*|^p)A^p.
\end{equation*}
\end{cor}

\vskip 5pt
\begin{proof} We apply Corollary \ref{corschurlast} to the pair of Hermitian operators in $\bM_{2n}$,
$$
X=\begin{pmatrix} 0&T^* \\ T&0\end{pmatrix}, \quad Y=\begin{pmatrix} 0&T \\ T^*&0\end{pmatrix}
$$
with  $A\oplus A$ in $\bM_{2n}^+$. We then obtain
$$
\left| \begin{pmatrix} 0&A(T\circ T^*)A \\ A(T\circ T^*)A&0\end{pmatrix} \right|^p \prec_{w\!\log} 
\begin{pmatrix} A(|T|^p\circ |T^*|^p)A&0 \\ 0&A(|T|^p\circ |T^*|^p)A\end{pmatrix}
$$
which is equivalent to the statement of our corollary. 
\end{proof}

\section{H\"older type inequalities}

\vskip 5pt Now we turn to H\"older's type inequalities. Fixing $p=1$ in Theorem \ref{th1}, we have the following special case.

\vskip 5pt
\begin{cor}\label{corspe} Let $A,B\in\bM_n^+$ and $Z\in\bM_n$. Then, for all symmetric norms and   $\alpha>0$, the map
$$
t \mapsto \left\| \left|A^{t}ZB^{t}\right|^{\alpha } \right\| 
$$
is  log-convex on $(-\infty,\infty)$.
\end{cor}

\vskip 5pt
This implies a fundamental fact, the L\"owner-Heinz inequality stating the operator monotonicity of $t^p$, $p\in(0,1)$.

\vskip 5pt
\begin{cor}  Let $A,B\in\bM_n^+$. If $A\ge B$, then $A^p\ge B^p$ for all $p\in(0,1)$.
\end{cor}

\vskip 5pt
\begin{proof} Corollary \ref{corspe} for the operator norm, with $Z=I$, $\alpha=2$, and the pair $A^{-1/2},B^{1/2}$ in place of the pair $A,B$ shows that
$
f(t)= \| A^{-t/2}B^{t}A^{-t/2}\|_{\infty}
$
is log-convex. Hence for $p\in(0,1)$, we have $f(p)\le f(1)^p f(0)^{1-p}$. Since $f(0)=1$ and by assumption $f(1)\le 1$, we obtain $f(p)\le 1$ and so $A^p\ge B^p$.
\end{proof}

Corollary \ref{corspe} entails a  H\"older inequality with a parameter. This inequality was first proved by Kosaki \cite[Theorem 3]{kosaki}. Here, we state it without the weight $Z$.

\vskip 5pt
\begin{cor}  Let $X,Y\in\bM_n$ and $p,q\ge 1$ such that $p^{-1}+q^{-1}=1$. Then, for all symmetric norms and   $\alpha>0$,  
$$
\left\| |XY|^{\alpha} \right\| \le  \left\|  |X|^{\alpha p}  \right\|^{1/p}  \left\|  |Y|^{q \alpha }  \right\|^{1/q}.
$$
\end{cor}

\vskip 5pt
\begin{proof} Let $A,B\in \bM_n^+$ with $B$ invertible. By replacing $B$ with $B^{-1}$
 and letting $Z=B$ in Corollary \ref{corspe} show that
$
t\mapsto  \| |A^tB^{1-t}|^{\alpha}\|
$
is log-convex on $(-\infty,\infty)$. Thus, for $t\in(0,1)$,
$
\| |A^tB^{1-t}|^{\alpha}\| \le \| A^{\alpha}\|^t\| B^{\alpha}\|^{1-t}
$.
Then, choose $A=|X|^{p}$, $B=|Y^*|^q$, $t=1/p$.
\end{proof}

\vskip 5pt More original H\"older's type inequalities are given in the next series of corollaries. 

\vskip 5pt
\begin{cor}\label{corcong} Let $A\in \bM_n^+$ and $Z\in\bM_{n,m}$. Then, for all symmetric norms and  $\alpha>0$, the map
$$
(p,t) \mapsto \left\| \left(Z^*A^{t/p}Z\right)^{\alpha p} \right\| 
$$
is jointly log-convex on $(0,\infty)\times(-\infty,\infty)$.
\end{cor}

\vskip 5pt
\begin{proof} By completing, if necessary, our matrices with some 0-entries, we may suppose $m=n$ and then apply Theorem \ref{th1} with $B=I$. \end{proof}

\vskip 5pt
\begin{cor}\label{corlittlewood} Let $a=(a_1,\cdots,a_m)$ and $w=(w_1,\cdots,w_m)$ be two $m$-tuples in $\bR^+$  and define, for all  $p>0$, 
$\| a\|_p:=(\sum_{i=1}^mw_ia_i^p)^{1/p}$. Then, for all $p,q>0$ and  $\theta\in(0,1)$, 
$$
\| a\|_{\frac{1}{\theta p +(1-\theta)q}} \le \| a\|_{\frac{1}{p}}^{\theta} \| a\|_{\frac{1}{q}}^{1-\theta}.
$$
\end{cor}

\vskip 5pt
\begin{proof} Fix $t=1$ and pick $A=\diag(a_1,\ldots a_m)$ and $Z^*=(w_1^{1/2},\ldots,w_m^{1/2})$ in the previous corollary. \end{proof}

\vskip 5pt
Corollary \ref{corlittlewood} is the classical log-convexity of $p\to \|\cdot\|_{1/p}$, or Littlewood's version of H\"older's inequality \cite[Theorem 5.5.1]{garling}. The next two corollaries, seemingly stronger but actually equivalent to Corollary \ref{corcong}, are also generalizations of this inequality.

\vskip 5pt
\begin{cor} Let $A_i\in\bM_n^+$ and $Z_i\in\bM_{n,m}$, $i=1,\ldots,k$. Then, for all symmetric norms and  $\alpha>0$, the map
$$
(p,t) \mapsto \left\| \left\{\sum_{i=1}^k Z_i^*A_i^{t/p}Z_i\right\}^{\alpha p} \right\| 
$$
is jointly log-convex on $(0,\infty)\times(-\infty,\infty)$.
\end{cor}

\vskip 5pt The unweighted case, $Z_i=I$ for all $i=1,\ldots,k$, is especially interesting. With $t=\alpha=1$, it is a matrix version of the unweighted Littlewood inequality.

\vskip 5pt
\begin{proof} Apply Corollary \ref{corcong} with $A=A_1\oplus\cdots\oplus A_k$ and $Z^*=(Z_1^*,\ldots,Z_k^*)$. \end{proof}

\vskip 5pt
\begin{cor}\label{poslin} Let $A\in \bM_m^+$ and let $\Phi:\bM_m^+\to\bM_n^+$ be a positive linear map. Then, for all symmetric norms and  $\alpha>0$, the map
$$
(p,t) \mapsto \left\| \left\{\Phi(A^{t/p})\right\}^{\alpha p} \right\| 
$$
is jointly log-convex on $(0,\infty)\times(-\infty,\infty)$.
\end{cor}

\vskip 5pt
\begin{proof} When restricted to the $*$-commutative subalgebra spanned by $A$, the map $\Phi$ has the form
\begin{equation}\label{decomp}\Phi(X)=\sum_{i=1}^m\sum_{j=1}^n Z_{i,j}^* XZ_{i,j}\end{equation}
for some rank 1 or 0 matrices $Z_{i,j}\in\bM_{m,n}$, $i=1,\ldots,m$, $j=1,\ldots,n$. So we are in the range of the previous corollary. To check the decomposition \eqref{decomp}, write the spectral decomposition
$A=\sum_{i=1}^m \lambda_i(A) E_i$
with rank one projections $E_i=x_i x_i^*$  for some column vectors $x_i\in \bM_{m,1}$ and set 
$Z_{i,j}= x_i R_{i,j}$ where $R_{i,j}\in\bM_{1,n}$ is the $j$-th row of $\Phi(E_i)^{1/2}$. \end{proof}

\vskip 5pt
The above proof shows a classical fact, a positive linear map on a commutative domain is completely positive. Our proof seems shorter than the ones in the literature.
We close this section with an application to Schur products.

\vskip 5pt
\begin{cor} Let $A,B\in \bM_n^+$. If $p\ge r\ge s\ge q$ and $p+q=r+s$, then, for all symmetric norms and   $\alpha>0$,
$$
\| \{A^r\circ B^s\}^{\alpha}\| \| \{A^s\circ B^r\}^{\alpha}\| \le \| \{A^p\circ B^q\}^{\alpha}\| \| \{A^q\circ B^p\}^{\alpha}\|
$$
and
$$
\| \{A^r\circ B^s\}^{\alpha}\| + \| \{A^s\circ B^r\}^{\alpha}\| \le \| \{A^p\circ B^q\}^{\alpha}\| +\| \{A^q\circ B^p\}^{\alpha}\|.
$$
\end{cor}

\vskip 5pt
\begin{proof} By a limit argument we may assume invertibility of $A$ and $B$. Let $w:=(p+q)/2$. We will show that the maps
\begin{equation}\label{eqshow}
t\mapsto \| \{A^{w+t} \circ B^{w-t}\}^{\alpha} \|, \quad t\mapsto \| \{A^{w-t} \circ B^{w+t}\}^{\alpha} \|
\end{equation}
are log-convex on $(-\infty,\infty)$. This implies that the functions
$$
f(t)=  \| \{A^{w+t} \circ B^{w-t}\}^{\alpha} \| \| \{A^{w-t} \circ B^{w+t}\}^{\alpha} \|
$$
and
$$
g(t)=  \| \{A^{w+t} \circ B^{w-t}\}^{\alpha} \| +\| \{A^{w-t} \circ B^{w+t}\}^{\alpha} \|
$$
are convex and even, hence nondecreasing on $[0,\infty)$. So we have
$$
f((r-s)/2)\le f((p-q)/2) \ {\text{and}} \ g((r-s)/2)\le g((p-q)/2)
$$
which prove the corollary.

To check the log-convexity of the maps \eqref{eqshow}, we  see the Schur product as a positive linear map acting on a tensor product,
$A\circ B= \Psi(A\otimes B)$. 
 By Corollary \ref{poslin}, the map
$$
t\mapsto \| \{\Phi(Z^t)\}^{\alpha} \|
$$
is log-convex on $(-\infty,\infty)$ for any positive matrix
 $Z\in\bM_n\otimes\bM_n$ and any positive linear map $\Phi:\bM_n\otimes\bM_n\to \bM_n$. Taking
$
Z=A\otimes B^{-1}
$
and
$$
\Phi(X)= \Psi(A^{w/2}\otimes B^{w/2}\cdot  X \cdot A^{w/2}\otimes B^{w/2})
$$
we obtain the log-convexity of the first map
$
t\mapsto \| \{A^{w+t}\circ B^{w-t}\}^{\alpha} \|
$
in $\eqref{eqshow}$. The log-convexity of the second one is similar.
\end{proof}

\section{Proof of Theorem \ref{th1}}

In the proof of the theorem, we will denote the $k$-th antisymmetric power $\wedge^k T$ of an operator $T$ simply as $T_k$. The symbol $\|\cdot\|_{\infty}$ stands for the usual operator norm while $\rho(\cdot)$ denotes the spectral radius. Given $A\in\bM_n^+$ 
we denote by $A^{\downarrow}$  the diagonal matrix with the eigenvalues of $A$ in decreasing order down to the diagonal, $A^{\downarrow}={\diag(\lambda_j(A)})$.

\begin{proof} Recall that if $A\in\bM_n^+$ is not invertible and $t\ge 0$, we define $A^{-t}$ as the generalized inverse of $A^t$, i.e., $A^{-t}:=(A+F)^{-t}E$ where $F$ is the projection onto the nullspace of $A$ and $E$ is the range projection of $A$. With this convention, replacing if necessary $Z$ by $EZE'$ where $E$ is the range projection of $A$ and $E'$ that  of $B$, we may and do assume that $A$ and $B$ are invertible.

Let
$$
g_k(t):= \prod_{j=1}^k \lambda_j(|A^{t}ZB^{t}|)= \|  A_k^{t}Z_kB_k^{t}\|_{\infty}
$$
Then
\begin{align*}g_k((t+s)/2) &= \| A_k^{(t+s)/2}Z_kB_k^{t+s}Z^*_kA_k^{(t+s)/2} \|_{\infty}^{1/2} \\
&=\rho^{1/2}( A_k^{t}Z_kB_k^{t+s}Z^*_kA_k^{s}) \\
&\le  \| A_k^{t}Z_kB_k^{t+s}Z^*_kA_k^{s} \|_{\infty}^{1/2} \\
&\le  \| A_k^{t}Z_kB_k^{t}\|_{\infty}^{1/2} \| B_k^{s}Z^*_kA_k^{s} \|_{\infty}^{1/2} \\
&=\{g_k(t)g_k(s)\}^{1/2}.
\end{align*}
Thus $t\mapsto g_k(t)$ is log-convex on $(-\infty,\infty)$ and so  $(p,t)\mapsto g_k^p(t/p)$ is jointly log-convex on $(0,\infty)\times(-\infty,\infty)$. Indeed, its logarithm $p\log g_k(t/p)$ 
is the perspective of the convex function $\log g_k(t)$, and hence is jointly convex. Therefore
\begin{equation}\label{equa1}
g_k^{(p+q)/2}\left(\frac{(t+s)/2}{(p+q)/2}\right) \le \{ g_k^p(t/p)g_k^q(s/q)\}^{1/2}
\end{equation}
for $k=1,2,\ldots,n$, with equality for $k=n$ as it then involves the determinant. This is equivalent to the log-majorization
$$
\left|A^{\frac{t+s}{p+q}}ZB^{\frac{t+s}{p+q}}\right|^{\frac{p+q}{2}} \prec_{\log}
|A^{t/p}ZB^{t/p}|^{\frac{p}{2}\downarrow}\left|A^{s/q}ZB^{s/q}\right|^{\frac{q}{2}\downarrow}
$$
which is equivalent, for any  $\alpha>0$, to the log-majorization
$$
\left|A^{\frac{t+s}{p+q}}ZB^{\frac{t+s}{p+q}}\right|^{\alpha \frac{p+q}{2}} \prec_{\log}
\left|A^{t/p}ZB^{t/p}\right|^{\frac{\alpha p}{2}\downarrow}\left|A^{s/q}ZB^{s/q}\right|^{\frac{\alpha q}{2}\downarrow}
$$
ensuring that
\begin{equation}\label{equa2}
\left\|  \left|A^{\frac{t+s}{p+q}}ZB^{\frac{t+s}{p+q}}\right|^{\alpha \frac{p+q}{2}}\right\|
\le \left\| \left|A^{t/p}ZB^{t/p}\right|^{\frac{\alpha p}{2}\downarrow}\left|A^{s/q}ZB^{s/q}\right|^{\frac{\alpha q}{2}\downarrow} \right\|
\end{equation}
for all symmetric norms. Thanks to the Cauchy-Schwarz inequality for symmetric norms, we then have
\begin{equation}\label{equa3}
\left\|  \left|A^{\frac{t+s}{p+q}}ZB^{\frac{t+s}{p+q}}\right|^{\alpha \frac{p+q}{2}}\right\|
\le \left\| \left|A^{t/p}ZB^{t/p}\right|^{\alpha p} \right\|^{1/2}\left\| \left|A^{t/q}ZB^{t/q}\right|^{\alpha q}\right\|^{1/2}
\end{equation}
which means that 
$$
(p,t) \mapsto \left\| \left|A^{t/p}ZB^{t/p}\right|^{\alpha p} \right\| 
$$
is jointly log-convex on $(0,\infty)\times(-\infty,\infty)$.
\end{proof}

\vskip 5pt 
Denote by $I_k$ the identity of $\bM_k$ and by $\det_k$
the determinant on $\bM_k$. If $n\ge k$, $\Theta(k,n)$ stands for the set of $n\times k$ isometry matrices $T$, i.e, $T^*T=I_k$. One easily checks the variational formula, for $A\in\bM_n$ and $k=1,\ldots,n$,
$$
\prod_{j=1}^k \lambda_j(|A|)=\max_{V,W\in\Theta(k,n)}\left| {\mathrm{det}}_k\, V^*AW\right|.
$$
From this formula follow two facts, Horn's inequality,
$$
\prod_{j=1}^k \lambda_j(|AB|) \le \prod_{j=1}^k \lambda_j(|A|)\lambda_j(|B|)
$$
for all $A,B\in\bM_n$ and $k=1,\dots,n$, and, making use of Schur's triangularization,  the inequality
$$
\prod_{j=1}^k \lambda_j(|AB|) \le \prod_{j=1}^k \lambda_j(|BA|)
$$
whenever $AB$ is normal (indeed, by Schur's theorem we may assume that $BA$ is upper triangular with the eigenvalues of $BA$, hence of $AB$, down to the diagonal and our variational formula then gives the above log-majorization). This shows that the proof of Theorem \ref{th1} can be written without the machinery of antisymmetric tensors.

The  novelty of this proof consists in using the perspective of a one variable convex function. One more perspective  yields the following variation of Corollary \ref{cornew}.

\vskip 5pt
\begin{cor} Let $A,B\in\bM_m^{+}$, let $Z\in\bM_m$. Then, for all symmetric norms and   $\alpha>0$, the map
$$
p \mapsto \left\| \left|A^{1/p}ZB^{1/p}\right|^{\alpha } \right\|^p 
$$
is  log-convex on $(0,\infty)$. 
\end{cor}

\vskip 5pt
\begin{proof} By Theorem \ref{th1} with fixed $p=1$, the map
$
t \mapsto \log\left( \| |A^{t}ZB^{t}|^{\alpha} \|\right)
$
is convex on $(0,\infty)$, thus its perpective
$$
(p,t) \mapsto p\log \left( \left\| |A^{t/p}ZB^{t/p}|^{\alpha} \right\|\right) = \log \left( \left\| |A^{t/p}ZB^{t/p}|^{\alpha} \right\|^p\right)
$$
is jointly log-convex on $(0,\infty)\times(0,\infty)$. Now fixing $t=1$ completes the proof.
\end{proof}

\vskip 5pt
From this corollary we may derive the next one exactly as Corollary \ref{poslin}
follows from Theorem \ref{th1}. This result is another noncommutative version of Littlewood's inequality (\cite[Theorem 5.5.1]{garling}).

\vskip 5pt
\begin{cor}  Let $\Phi:\bM_m\to \bM_n$ be a positive linear map and let $A\in\bM_m^+$. Then, for all symmetric norms and  $\alpha>0$,
 the map
$$
p \mapsto \left\| \Phi^{\alpha }(A^{1/p}) \right\|^p
$$
is  log-convex on $(0,\infty)$. 
\end{cor}

\section{Hilbert space operators}

\vskip 5pt
In this section we give a version  
 of Theorem \ref{th1}  for the algebra $\bB$ of bounded linear operators on a separable, infinite dimensional 
Hilbert space ${\mathcal{H}}$. 
We first include a brief treatment of symmetric norms for operators in $\bB$. Our approach does not require to discuss any underlying ideal, we refer the reader to \cite[Chapter 2]{Simon} for a much more complete discussion.

We may  define symmetric norms on $\bB$ in a closely related way to the finite dimensional case as follows. Let  $\bF$ be the set of finite rank operators and $\bF^+$ its positive part.

\vskip 5pt\noindent
\begin{definition}\label{def}A symmetric norm $\|\cdot\|$ on $\bB$ is a functional taking value in $[0,\infty]$ such that:
\begin{itemize}
\item[(1)] $\|\cdot \|$ induces a norm on $\bF$.
\item[(2)] If $\{X_n\}$ is a  sequence  in $\bF^+$ strongly increasing to $X$, then $\|X\|=\lim_n \| X_n\|$.
\item[(3)] $\| KZL\| \le \| Z\|$ for all $Z\in \bB$ and all contractions $K,L\in\bB$.
\end{itemize}
\end{definition}

\vskip 5pt
The reader familiar to the  theory of symmetrically normed ideals may note that our definition of a symmetric norm is equivalent to the usual one. More precisely, restricting $\|\cdot\|$ to the set where it takes finite values, Definition \ref{def} yields the classical notion of a symmetric norm defined on its maximal ideal.

Definition \ref{def} shows that a symmetric norm on $\bB$ induces a symmetric norm on $\bM_n$ for each $n$, say $\|\cdot\|_{\bM_n}$. In fact $\|\cdot\|$ can be regarded as a limit of the norms $\|\cdot\|_{\bM_n}$, see Lemma \ref{lem1} for a precise statement, so that basic properties of symmetric norms on $\bM_n$ can be extended to symmetric norms on $\bB$. For instance the Cauchy-Schwarz inequality also holds   for symmetric norms on $\bB$, (with possibly the $\infty$ value) as well as the Ky Fan principle for $A,B\in\bB^+$:
If $A\prec_w B$, then $\| A\| \le \| B\|$ for all symmetric norms. In fact, even for a noncompact operator $A\in\bB^+$, the sequence  $\{\lambda_j(A)\}_{j=1}^{\infty}$  and the corresponding diagonal operator $A^{\downarrow}={\diag(\lambda_j(A)})$  are well defined, via the minmax formulae (see \cite[Proposition 1.4]{Hiai})
$$
\lambda_j(A) = \inf_E\{ \| EAE \|_{\infty} \ : E  \ {\mathrm{ projection\ with\ }}{\mathrm{rank}}(I-E)=j-1 \}
$$
 The  Ky Fan principle then still holds for $A,B\in\bB^+$ by Lemma \ref{lem1} and the obvious property
$$
\lambda_j(A)=\lim_{n\to\infty} \lambda_j(E_nAE_n)
$$
for all sequences of finite rank projections $\{E_n\}_{n=1}^{\infty}$ strongly converging to the identity. 
Note also that we still have $\| \wedge^k A\|_{\infty} =\prod_{j=1}^k \lambda_j(A)$.

Thus we have the same tools as in the matrix case and we will be able to adapt the proof of Theorem \ref{th1}  for  $\bB$. The infinite dimensional version of Theorem \ref{th1}
is the following statement.

\vskip 5pt
\begin{theorem}\label{th2} Let $A,B\in\bB^{+}$, let $Z\in\bB$. Then, for all symmetric norms and   $\alpha>0$, the map
$$
(p,t) \mapsto \left\| \left|A^{t/p}ZB^{t/p}\right|^{\alpha p} \right\| 
$$
is jointly log-convex on $(0,\infty)\times(0,\infty)$. This map takes its finite values in the  open quarter-plan 
$$\Omega(p_0,t_0)=\{ (p,t) \ | \ p>p_0, \ t>t_0\}$$
for some $p_0,t_0\in [0,\infty]$, or on its  closure $\overline{\Omega}(p_0,t_0)$.
\end{theorem}

\vskip 5pt
Note that, contrarily to Theorem \ref{th1}, we confine the variable  $t$ to the positive half-line. Indeed, when dealing with a symmetric norm, the operators $A$ and $B$ are often compact, so that, for domain reasons,  we cannot consider two unbounded operators such as $A^{-1}$ and $B^{-1}$.

\vskip 5pt
\begin{proof} Note that $AZB=0$ if and only if $A^qZB^q=0$, for any $q>0$. In this case, our map is the 0-map, and  its logarithm with constant value $-\infty$ can be regarded as convex. Excluding this trivial case, our map takes values in $(0,\infty]$ and it makes sense to consider the log-convexity property. 
   We may reproduce the  proof of Theorem \ref{th1} and obtain \eqref{equa1}
for all $k=1,2,\cdots$. This leads to weak-logmajorizations and so to a weak majorization equivalent (Ky Fan's principle in $\bB$) to \eqref{equa2}, with possibly the $\infty$ value on the right side or both sides. The Cauchy-Schwarz inequality for symmetric norms in $\bB$ yields \eqref{equa3} (possibly with the $\infty$ value). Therefore our map is jointly log-convex. To show that the domain where it takes finite values is $\Omega(p_0,t_0)$ or $\overline{\Omega}(p_0,t_0)$, it suffices to show the following two implications: 

\vskip 5pt
{\it Let $0<t<s$ and $0<p<q$. If $\left\| \left|A^{t/p}ZB^{t/p}\right|^{\alpha p} \right\| <\infty$, then
\begin{itemize}
\item[{\rm(i)}] $\left\| \left|A^{s/p}ZB^{s/p}\right|^{\alpha p} \right\| <\infty$, and
\item[{\rm (ii)}] $\left\| \left|A^{t/q}ZB^{t/q}\right|^{\alpha q} \right\| <\infty$.
\end{itemize} }

\vskip 5pt
Since $0<t<s$ ensures that, for some constant $c=c(s,t)>0$,
$$
\lambda_j(|A^{t/p}ZB^{t/p}|) \ge c\lambda_j(|A^{s/p}ZB^{s/p}|)
$$
for all $j=1,2,\ldots$, we obtain (i). To obtain (ii) we may assume that $Z$ is a contraction. Then arguing as in the proof of Corollary \ref{cor1} we see that the finite value map
$$p\mapsto  \left\| \left|A^{t/p}ZB^{t/p}\right|^{\alpha p} \right\| $$
is nonincreasing for all Ky-Fan norms. Thus this map is also nonincreasing for all symmetric norms. This gives (ii).
\end{proof}

\vskip 5pt
Exactly as in the matrix case, we can derive the following two corollaries.

\vskip 5pt
\begin{cor}  Let $A,B\in\bB^+$ and $p\ge 1$. Then,  for all contractions $Z\in \bB$,
$$
   (AZ^*BZA)^p \prec_{w\!\log}   A^{p}Z^*B^{p}ZA^{p}.
$$
\end{cor}

\vskip 5pt
\begin{cor}  Let $A,B\in\bB^+$ and  let $Z\in\bB$ be a contraction. Assume that at least one of these three operators is compact. Then, if $p\ge 1$ and 
$f(t)$ is e-convex and nondecreasing,
$$
  {\mathrm{Tr\,}} f((AZ^*BZA)^p) \le  {\mathrm{Tr\,}} f( A^{p}Z^*B^{p}ZA^{p}).
$$
\end{cor}

\vskip 5pt
Here, we use the fact that for $X\in\bK^+$ and a nondecreasing continuous function $f:[0,\infty)\to (-\infty,\infty)$, we can define ${\mathrm{Tr\,}} f(X)$ as an element in  $[-\infty,\infty]$ by 
$$
  {\mathrm{Tr\,}} f(X)=\lim_{k\to \infty} \sum_{j=1}^k f(\lambda_j(X)).
$$

\vskip 5pt
Given a symmetric norm $\|\cdot\|$ on $\bB$, the set where $\|\cdot\|$ takes a finite value is an ideal. We call it the maximal ideal of $\|\cdot\|$ or the domain of $\|\cdot\|$. From Theorem \ref{th2} we immediately infer our last corollary.

\vskip 5pt
\begin{cor}  Let $A,B\in\bB^+$ and   $Z\in\bB$. Suppose that $AZB\in\bJ$, the domain of a symmetric norm. Then, for all $q\in(0,1)$, we also have $|A^qZB^q|^{1/q}\in\bJ$. 
\end{cor}

\vskip 5pt
Following \cite[Chapter 2]{Simon}, we denote by $\bJ^{(0)}$ the $\|\cdot\|$-closure of the finite rank operators. In most cases $\bJ=\bJ^{(0)}$, however the strict inclusion  $\bJ^{(0)}\subset_{\neq} \bJ$ may happen. We do not know whether we can replace in the last corollary  $\bJ$ by $\bJ^{(0)}$.

We close our article with two simple lemmas  and show how the Cauchy-Schwarz inequality for the infinite dimensional case follows from the matrix case.

\vskip 5pt\noindent
\begin{lemma}\label{lem1} Let $\|\cdot \|$ be a symmetric norm on $\bB$ and let $\{E_n\}_{n=1}^{\infty}$ be an increasing sequence of finite rank projections in $\bB$, strongly converging to $I$. Then, for all $X\in \bB$, $\| X\|=\lim_n \| E_nXE_n\|$.
\end{lemma}

\vskip 5pt
\begin{proof} We first show that $\| E_nX\| \to \| X\|$ as $n\to \infty$.  Since  
$\| E_n X\|=\| (X^*E_nX)^{1/2}\|$ and $(X^*E_nX)^{1/2} \nearrow |X|$ by operator monotonicity of $t^{1/2}$, we obtain $\lim_n\| E_nX\|= \| X\|$ by Definition \ref{def}(2). Similarly, 
$\lim_k\| E_nXE_k\|=\| E_nX\|$, and so $\lim_n\| E_nXE_{k(n)}\|=\| X\|$, and thus,  by Definition \ref{def}(3), $\lim_p\| E_pXE_p\|=\| X\|$.
\end{proof}

\vskip 5pt\noindent
\begin{lemma}\label{lem2} Let $\|\cdot \|$ be a symmetric norm on $\bB$ and let $\{E_n\}_{n=1}^{\infty}$ and $\{F_n\}_{n=1}^{\infty}$ be two increasing sequences of finite rank projections in $\bB$, strongly converging to $I$. Then, for all $X\in \bB$, $\| X^*X\|=\lim_n \| E_nX^*F_nXE_n\|$.
\end{lemma}

\vskip 5pt
\begin{proof} By Definition \ref{def}(2)-(3), the map $n\mapsto \| E_nX^*F_nXE_n\|$ is nondecreasing. By Definition \ref{def}(2), for any integer $p$, its limit is greater than or equal $\| E_pX^*XE_p\|$. By Lemma \ref{lem1}, the limit is precisely $\| X^*X\|$.
\end{proof}

\vskip 5pt
Let $X,Y\in\bB$, let $\|\cdot\|$ be a symmetric norm on $\bB$, and let $\{E_n\}_{n=1}^{\infty}$ be as in the above lemma. Let $F_n$ be the range projection of $YE_n$. We have by Lemma \ref{lem1}
$$
\| X^*Y\| =\lim_n \| E_n X^*Y E_n \| =\lim_n  \| E_n X^*F_nY E_n \|.
$$
Let ${\mathcal{H}}_n$ be the sum of the ranges of $E_n$ and  $F_n$. This is a finite dimensional subspace, say $\dim{\mathcal{H}}_n=d(n)$. Applying the Cauchy-Schwarz inequality for a symmetric norm on $\bM_{d(n)}$, we obtain, thanks to Lemma \ref{lem2},
\begin{align*}
\| X^*Y\| &=\lim_n  \| E_n X^*F_nY E_n\|_{\bM_{d(n)}}  \\
&\le \lim_n  \| E_n X^*F_nX E_n\|_{\bM_{d(n)}}^{1/2}  \| E_n Y^*F_n YE_n\|_{\bM_{d(n)}}^{1/2}  \\
&= \|  X^*X \|^{1/2} \|  Y^*Y \|^{1/2}. 
\end{align*}
Thus the Cauchy-Schwarz inequality for a symmetric  norm on $\bB$ follows from the Cauchy-Schwarz inequality for symmetric norms on $\bM_{n}$. Of course, the two previous lemmas and this discussion are rather trivial, but we wanted to stress on the fact that Theorem \ref{th2} is essentially of finite dimensional nature. However, it would be also desirable to extend these results in the setting of a semifinite von Neumann algebra.

\vskip 5pt
\noindent
Laboratoire de math\'ematiques, 
Universit\'e de Franche-Comt\'e, 
25 000 Besan\c{c}on, France.
Email: jcbourin@univ-fcomte.fr

\vskip 5pt\noindent
 Department of mathematics, 
Kyungpook National University, 
 Daegu 702-701, Korea.
Email: eylee89@knu.ac.kr

\end{document}